\documentclass[12pt]{article}
\usepackage{amsfonts}
\usepackage{full page,
amssymb, amscd, amsmath,graphicx, 
}

 \title{{\bf Affine Lie algebras and tensor categories}}
 \author{Yi-Zhi Huang}
    \date{}
    \begin{document}
    \bibliographystyle{alpha}

\newtheorem{thm}{Theorem}[section]
\newtheorem{defn}[thm]{Definition}
\newtheorem{prop}[thm]{Proposition}
\newtheorem{cor}[thm]{Corollary}
\newtheorem{lemma}[thm]{Lemma}
\newtheorem{rema}[thm]{Remark}
\newtheorem{app}[thm]{Application}
\newtheorem{prob}[thm]{Problem}
\newtheorem{conv}[thm]{Convention}
\newtheorem{conj}[thm]{Conjecture}
\newtheorem{cond}[thm]{Condition}
    \newtheorem{exam}[thm]{Example}
\newtheorem{assum}[thm]{Assumption}
     \newtheorem{nota}[thm]{Notation}
\newcommand{\halmos}{\rule{1ex}{1.4ex}}
\newcommand{\pfbox}{\hspace*{\fill}\mbox{$\halmos$}}

\newcommand{\nn}{\nonumber \\}

 \newcommand{\res}{\mbox{\rm Res}}
 \newcommand{\ord}{\mbox{\rm ord}}
\renewcommand{\hom}{\mbox{\rm Hom}}
\newcommand{\edo}{\mbox{\rm End}\ }
 \newcommand{\pf}{{\it Proof.}\hspace{2ex}}
 \newcommand{\epf}{\hspace*{\fill}\mbox{$\halmos$}}
 \newcommand{\epfv}{\hspace*{\fill}\mbox{$\halmos$}\vspace{1em}}
 \newcommand{\epfe}{\hspace{2em}\halmos}
\newcommand{\nord}{\mbox{\scriptsize ${\circ\atop\circ}$}}
\newcommand{\wt}{\mbox{\rm wt}\ }
\newcommand{\swt}{\mbox{\rm {\scriptsize wt}}\ }
\newcommand{\lwt}{\mbox{\rm wt}^{L}\;}
\newcommand{\rwt}{\mbox{\rm wt}^{R}\;}
\newcommand{\slwt}{\mbox{\rm {\scriptsize wt}}^{L}\,}
\newcommand{\srwt}{\mbox{\rm {\scriptsize wt}}^{R}\,}
\newcommand{\clr}{\mbox{\rm clr}\ }
\newcommand{\tr}{\mbox{\rm Tr}}
\newcommand{\C}{\mathbb{C}}
\newcommand{\Z}{\mathbb{Z}}
\newcommand{\R}{\mathbb{R}}
\newcommand{\Q}{\mathbb{Q}}
\newcommand{\N}{\mathbb{N}}
\newcommand{\CN}{\mathcal{N}}
\newcommand{\F}{\mathcal{F}}
\newcommand{\I}{\mathcal{I}}
\newcommand{\V}{\mathcal{V}}
\newcommand{\one}{\mathbf{1}}
\newcommand{\BY}{\mathbb{Y}}
\newcommand{\ds}{\displaystyle}

        \newcommand{\ba}{\begin{array}}
        \newcommand{\ea}{\end{array}}
        \newcommand{\be}{\begin{equation}}
        \newcommand{\ee}{\end{equation}}
        \newcommand{\bea}{\begin{eqnarray}}
        \newcommand{\eea}{\end{eqnarray}}
         \newcommand{\lbar}{\bigg\vert}
        \newcommand{\p}{\partial}
        \newcommand{\dps}{\displaystyle}
        \newcommand{\bra}{\langle}
        \newcommand{\ket}{\rangle}

        \newcommand{\ob}{{\rm ob}\,}
        \renewcommand{\hom}{{\rm Hom}}

\newcommand{\A}{\mathcal{A}}
\newcommand{\Y}{\mathcal{Y}}

\newcommand{\dlt}[3]{#1 ^{-1}\delta \bigg( \frac{#2 #3 }{#1 }\bigg) }

\newcommand{\dlti}[3]{#1 \delta \bigg( \frac{#2 #3 }{#1 ^{-1}}\bigg) }

 \makeatletter
\newlength{\@pxlwd} \newlength{\@rulewd} \newlength{\@pxlht}
\catcode`.=\active \catcode`B=\active \catcode`:=\active
\catcode`|=\active
\def\sprite#1(#2,#3)[#4,#5]{
   \edef\@sprbox{\expandafter\@cdr\string#1\@nil @box}
   \expandafter\newsavebox\csname\@sprbox\endcsname
   \edef#1{\expandafter\usebox\csname\@sprbox\endcsname}
   \expandafter\setbox\csname\@sprbox\endcsname =\hbox\bgroup
   \vbox\bgroup
  \catcode`.=\active\catcode`B=\active\catcode`:=\active\catcode`|=\active
      \@pxlwd=#4 \divide\@pxlwd by #3 \@rulewd=\@pxlwd
      \@pxlht=#5 \divide\@pxlht by #2
      \def .{\hskip \@pxlwd \ignorespaces}
      \def B{\@ifnextchar B{\advance\@rulewd by \@pxlwd}{\vrule
         height \@pxlht width \@rulewd depth 0 pt \@rulewd=\@pxlwd}}
      \def :{\hbox\bgroup\vrule height \@pxlht width 0pt depth
0pt\ignorespaces}
      \def |{\vrule height \@pxlht width 0pt depth 0pt\egroup
         \prevdepth= -1000 pt}
   }
\def\endsprite{\egroup\egroup}
\catcode`.=12 \catcode`B=11 \catcode`:=12 \catcode`|=12\relax
\makeatother

\def\hboxtr{\FormOfHboxtr} 
\sprite{\FormOfHboxtr}(25,25)[0.5 em, 1.2 ex] 

:BBBBBBBBBBBBBBBBBBBBBBBBB | :BB......................B |
:B.B.....................B | :B..B....................B |
:B...B...................B | :B....B..................B |
:B.....B.................B | :B......B................B |
:B.......B...............B | :B........B..............B |
:B.........B.............B | :B..........B............B |
:B...........B...........B | :B............B..........B |
:B.............B.........B | :B..............B........B |
:B...............B.......B | :B................B......B |
:B.................B.....B | :B..................B....B |
:B...................B...B | :B....................B..B |
:B.....................B.B | :B......................BB |
:BBBBBBBBBBBBBBBBBBBBBBBBB |

\endsprite
\def\shboxtr{\FormOfShboxtr} 
\sprite{\FormOfShboxtr}(25,25)[0.3 em, 0.72 ex] 

:BBBBBBBBBBBBBBBBBBBBBBBBB | :BB......................B |
:B.B.....................B | :B..B....................B |
:B...B...................B | :B....B..................B |
:B.....B.................B | :B......B................B |
:B.......B...............B | :B........B..............B |
:B.........B.............B | :B..........B............B |
:B...........B...........B | :B............B..........B |
:B.............B.........B | :B..............B........B |
:B...............B.......B | :B................B......B |
:B.................B.....B | :B..................B....B |
:B...................B...B | :B....................B..B |
:B.....................B.B | :B......................BB |
:BBBBBBBBBBBBBBBBBBBBBBBBB |

\endsprite

\vspace{2em}



\date{}
\bibliographystyle{alpha}
\maketitle

\begin{abstract}
We review briefly the existing
vertex-operator-algebraic constructions of various tensor category
structures on module categories for affine Lie
algebras. We discuss the results that were first conjectured in the work of 
Moore and Seiberg and led us to the construction
of the modular tensor category structure in the positive integral level case. 
Then we review the existing constructions and results in the following three cases:
(i)  the level plus the dual Coxeter number is not a nonnegative rational number,
(ii)  the level is a positive integer and (iii)  the level 
is an admissible number. We also present several open problems.
\end{abstract}

\section{Introduction}

In 1988, Moore and Seiberg \cite{MS1} \cite{MS2} discovered that rational conformal field theories 
have properties analogous to the properties of tensor categories satisfying 
additional conditions. The results of Moore and Seiberg were obtained based on 
the operator product expansion and modular invariance of chiral vertex operators 
(intertwining operators). Mathematically, these two fundamental properties
had been deep and important conjectures for many years. But it was 
the Verlinde formula and modular tensor category structure in the 
case of affine Lie algebras at positive integral levels
(the two most well-known consequences of these two conjectures mentioned above)
that first attracted the attentions of mathematicians. 

To use the tensor category structures to study the representation theory of affine Lie algebras
and to solve related problems in algebra, topology, mathematical physics and other related 
areas, the first mathematical problem is to give constructions of these tensor categories. 
Kazhdan and Lusztig \cite{KL1}--\cite{KL5}
first constructed a rigid braided tensor category 
structure on a category $\mathcal{O}_{\ell}$ of modules for the affine Lie algebra 
$\hat{\mathfrak{g}}$ of a finite-dimensional complex 
simple Lie algebra $\mathfrak{g}$ in the case 
that the level $\ell$ plus the dual Coxeter number is not a nonnegative rational number. 
Then several works, including those of Beilinson, Mazur and Feigin \cite{BFM} and
Lepowsky and the author \cite{HL-affine},  constructed 
a braided tensor category structure using different approaches 
on a category $\widetilde{O}_{\ell}$ 
of $\hat{\mathfrak{g}}$-modules in the case that the level $\ell$ is a positive integer. 
The modular tensor category structure on $\widetilde{O}_{\ell}$ 
in the case that $\ell$ is a positive integer
was constructed by the author \cite{H-rigidity} as a special case 
of a general construction of the modular tensor category structures
conjectured by Moore and Seiberg \cite{MS2} for vertex operator algebras corresponding 
to rational conformal field theories. After a gap was filled later in \cite{F3}
using the Verlinde formula
proved by Faltings \cite{Fa}, Teleman \cite{Te} and the author \cite{H-Verlinde}, 
Finkelberg's work \cite{F1}--\cite{F3} 
can also be reinterpreted as giving a construction of the modular tensor category 
structure on $\widetilde{O}_{\ell}$  except for a few cases, including the important case
of $\mathfrak{g}=\mathfrak{e}_{8}$ and $\ell=2$. 
In the case that the level $\ell$ is admissible, 
Creutzig, Yang and the author \cite{CHY} constructed a braided tensor 
category structure on a semisimple category $\widetilde{O}_{\ell, {\rm \scriptsize ord}}$ 
of $\hat{\mathfrak{g}}$-modules. In the special case
of $\mathfrak{g}=\mathfrak{sl}_{2}$, the rigidity 
of $\widetilde{O}_{\ell, {\rm \scriptsize ord}}$ was proved and the problem 
of  determining when $\widetilde{O}_{\ell, {\rm \scriptsize ord}}$ is a modular tensor category 
was solved completely in \cite{CHY}. Recently, Creutzig proved in \cite{C} 
that when $\mathfrak{g}$
is simple-laced and the level $\ell$ is admissible, 
the braided tensor category structure on $\widetilde{O}_{\ell, {\rm \scriptsize ord}}$ 
is rigid and thus $\widetilde{O}_{\ell, {\rm \scriptsize ord}}$ is a ribbon tensor
category.

These constructions are all implicitly or explicitly 
vertex-operator-algebraic. In the case that $\mathfrak{g}=\mathfrak{sl}_{n}$ and 
the level is a positive integer, 
Kawahigashi, Longo and M\"{u}ger \cite{KLM} 
gave a construction of the modular tensor category structure on $\widetilde{O}_{\ell}$ 
using the approach of conformal nets.

In all the cases discussed above, 
not only have the braided tensor category structures been constructed,
but also the underlying vertex tensor category structures or modular vertex tensor category structures
(see \cite{HL1}--\cite{HL5}, \cite{tensor4}, \cite{HL-affine},
\cite{HLZ1}--\cite{HLZ9}, \cite{H-rigidity}, \cite{Z}, \cite{H-Z-correction} and \cite{CHY}). 
We would like to emphasize that from a vertex tensor category, we obtain naturally
a braided tensor category (see \cite{HL2} and \cite{HL6}) but in general, 
a braided tensor category alone is not enough to construct a vertex tensor category.
This fact corresponds to the fact that from a two-dimensional conformal
field theory, we obtain naturally a three-dimensional topological field theory
but in general,  a three-dimensional topological field theory alone is not 
enough to construct a two-dimensional conformal field theory. 
Vertex tensor categories (see \cite{HL2}), rigid vertex tensor categories (see \cite{H-rigidity}) 
or modular vertex tensor categories (see \cite{H-rigidity}) contain 
much more information than braided tensor categories, ribbon tensor categories
or modular tensor categories, respectively.

In this paper, we review briefly the existing vertex-operator-algebraic
constructions of tensor category structures
on module categories for affine Lie algebras. We also present several  
open problems in this paper. For a review on the approach of 
conformal nets, including in particular
the work \cite{KLM}, we refer the reader to
\cite{K}. 

We recall some basic definitions and constructions in the representation theory of 
affine Lie algebras in the next section. In Section 3, we 
review briefly the important work \cite{MS1} 
\cite{MS2} of Moore and Seiberg. The review of the existing vertex-operator-algebraic
constructions of various tensor category
structures on module categories for affine Lie
algebras are given in Section 4. In Subsection 4.1, we recall some basic terms 
in the theory of tensor categories. In Subsections 4.2, 4.3 and 4.4, we review the 
existing constructions in the case that the level plus the dual Coxeter number 
is not a nonnegative rational number, in
the case that the level is a positive integer and  in the case that the level 
is an admissible number, respectively. The open problems in these cases are also 
given in these subsections. We state the open problem for the remaining case
in Subsection 4.5.

\paragraph{Acknowledgments}
This paper is an expanded and revised version of the slides of the  author's 
talk at the ``10th Seminar on Conformal Field Theory: A conference on 
Vertex Algebras and Related Topics'' held at the Research Institute for Mathematical Sciences,
Kyoto University, Kyoto, April 23--27, 2018. The author is grateful
to the organizers Tomoyuki Arakawa, Peter Fiebig,
Nils Scheithauer, Katrin Wendland and Hiroshi Yamauchi for their invitation
and to Nils Scheithauer for the financial support.  The author would also like to thank Thomas
Creutzig for comments.

\section{Affine Lie algebras and modules}

We recall the basic definitions and constructions in the representation theory of 
affine Lie algebras. 

Let $\mathfrak{g}$ be  a finite-dimensional complex simple Lie algebra of rank $r$ 
and  $(\cdot, \cdot)$ the invariant symmetric bilinear form on $\mathfrak{g}$.
The affine Lie algebra $\hat{\mathfrak{g}}$ is
the vector space $\mathfrak{g} \otimes \C[t, t^{-1}] \oplus \C{\bf k}$ equipped with the bracket operation
\begin{align*}
[a \otimes t^m, b \otimes t^n] &= [a, b]\otimes t^{m+n} + (a, b)m\delta_{m+n,0}{\bf k},\\
[a \otimes t^m, {\bf k}] &= 0,
\end{align*}
for $a, b \in \mathfrak{g}$ and $m, n \in \Z$. 
Let $\hat{\mathfrak{g}}_{\pm} = \mathfrak{g}\otimes t^{\pm 1}\C[t^{\pm 1}]$. Then
$$\hat{\mathfrak{g}} = \hat{\mathfrak{g}}_{-} \oplus \mathfrak{g} \oplus \C{\bf k} \oplus \hat{\mathfrak{g}}_{+}.$$

If ${\bf k}$ acts as a complex number $\ell$ on a $\hat{\mathfrak{g}}$-module, 
then $\ell$ is called the level of this module.

Let $h$ and $h^{\vee}$ be the Coxeter number and dual Coxeter number, respectively,
of $\mathfrak{g}$. 
Let $M$ be a $\mathfrak{g}$-module and let $\ell \in \C$.
Assume that M can be decomposed as a direct sum of 
generalized eigenspaces of the Casimir operator. This decomposition gives a $\C$-grading to $M$.
In the case that $\ell+h^{\vee}\ne 0$, by using this decomposition and 
defining the (conformal) weight of the generalized eigenspace 
of the Casimir operator with eigenvalue $n$ to be $\frac{n}{2(\ell+h^{\vee})}$, 
we obtain a $\C$-grading of $M$ by (conformal) weights. Let $\hat{\mathfrak{g}}_{+}$ act on $M$ trivially and let
$\mathbf{k}$ act as the scalar multiplication by $\ell$. 
Then $M$ becomes a $\mathfrak{g} \oplus \C{\bf k} \oplus \hat{\mathfrak{g}}_{+}$-module
and we have a $\C$-graded induced $\hat{\mathfrak{g}}$-module
$$\widehat{M}_{\ell} = U(\hat{\mathfrak{g}})
\otimes_{U(\mathfrak{g}\oplus \C\mathbf{k}\oplus \hat{\mathfrak{g}}_{+})}M,$$
where the $\C$-grading is given by the $\C$-grading on $M$ and 
the grading on  $U(\hat{\mathfrak{g}})$ induced from 
the grading on $\hat{\mathfrak{g}}$.

Let $\mathfrak{h}$ be a Cartan subalgebra of $\mathfrak{g}$. 
For $\lambda\in \mathfrak{h}^{*}$,
let $L(\lambda)$ be the irreducible highest weight $\mathfrak{g}$-module 
with the highest weight $\lambda$. We use $M(\ell, \lambda)$ to denote the 
$\hat{\mathfrak{g}}$-module $\widehat{L(\lambda)}_{\ell}$. 
Let $J(\ell, \lambda)$ be the maximal proper submodule of $M(\ell, \lambda)$ 
and $L(\ell, \lambda) = M(\ell, \lambda)/J(\ell, \lambda)$. 
Then $L(\ell, \lambda)$ is the unique irreducible graded $\hat{\mathfrak{g}}$-module 
such that ${\bf k}$ acts as $\ell$ and the space of all elements annihilated by 
$\hat{\mathfrak{g}}_{+}$ is isomorphic to the $\mathfrak{g}$-module $L(\lambda)$.

Frenkel and Zhu  \cite{FZ} gave both $M(\ell, 0)$ and $L(\ell, 0)$ natural structures 
of  vertex operator algebras. Moreover, they gave in \cite{FZ}
both $M(\ell, \lambda)$ and $L(\ell, \lambda)$ structures of $M(\ell, 0)$-modules and  
$L(\ell, \lambda)$ a structure of an $L(\ell, 0)$-module 
for dominant integral $\lambda$.

For $\ell \in \C$ such that $\ell+h^{\vee}\not\in \Q_{\ge 0}$, let $\mathcal{O}_{\ell}$
be the category of all the $\hat{\mathfrak{g}}$-modules of level $\ell$ having a finite
composition series all of whose irreducible subquotients are of the form
$L(\ell, \lambda)$ for dominant integral $\lambda\in \mathfrak{h}^{*}$.

For $\ell\in \Z_{+}$, let $\widetilde{\mathcal{O}}_{\ell}$ be the category of 
$\hat{\mathfrak{g}}$-modules of level $\ell$ that are isomorphic to direct sums of 
irreducible $\hat{\mathfrak{g}}$-modules of the form $L(\ell, \lambda)$
for dominant integral $\lambda\in \mathfrak{h}^{*}$ 
such that $(\lambda, \theta)\le \ell$, where $\theta$ is the highest root of $\mathfrak{g}$.

Admissible modules for affine Lie algebras were studied first by Kac and Wakimoto
\cite{KW1} \cite{KW2}. 
The level of these modules are called admissible numbers.
Let $\ell$ be an admissible number, that is,
$\ell+h^{\vee}=\frac{p}{q}$ for $p, q\in \Z_{+}$, $(p, q)=1$, 
$p\ge h^{\vee}$ if $(r^{\vee}, q)=1$ and $p\ge h$ if $(r^{\vee}, q)=r^{\vee}$, where
$r^{\vee}$ is the "lacety" or 
lacing number of $\mathfrak{g}$,
that is, the maximum number of edges in the Dynkin diagram of $\mathfrak{g}$. 
Let $\mathcal{O}_{\ell, {\rm \scriptsize ord}}$ be the 
category of $\hat{\mathfrak{g}}$-modules of level $\ell$ that are isomorphic to direct sums of 
irreducible modules for the 
vertex operator algebra $L(\ell, 0)$.

\section{Operator product expansion, modular invariance and Moore-Seiberg equations}

The important work \cite{MS1}
\cite{MS2} of Moore-Seiberg on two-dimensional conformal field theory 
 in 1988 led to a conjecture: The category $\widetilde{O}_{\ell}$ for $\ell\in \Z_{+}$
has a natural structure of a modular tensor category in the sense of Turaev \cite{T1}
\cite{T2}. This conjecture was proved by the author in \cite{H-rigidity} in 2005 as a special case 
of a general construction of modular tensor categories for vertex operator algebras 
satisfying natural existence-of-nondegenerate-invariant-bilinear form, 
positive-energy, finiteness and reductivity conditions.

Since there are still a lot of misunderstandings about this work of Moore and Seiberg, here
we clarify what were the main assumptions in \cite{MS1} and \cite{MS2} and what were 
in fact proved in \cite{MS1} and \cite{MS2}. In \cite{MS2}, Moore and Seiberg formulated two 
major conjectures
on chiral rational conformal field theories:  The operator product expansion
of chiral vertex operators (see the end of Page 190 in \cite{MS2}) 
and the modular invariance of chiral vertex operators (see the beginning 
of the second paragraph on Page 200 in \cite{MS2}). In \cite{MS2}, no 
proofs or explanations or discussions
were given on why the operator product expansion and modular invariance of 
the chiral vertex operators formulated in \cite{MS2} are true. These were 
the most fundamental hypotheses, not results, of this work of Moore and Seiberg. 
Mathematically, chiral vertex operators are called intertwining operators in \cite{FHL}. 
Conformal field theory can be viewed as the study of intertwining operators. 
The two major hypotheses above are clearly major mathematical conjectures
which turned out to be more difficult to prove than the consequences derived from them in
the work \cite{MS1} and \cite{MS2}.

From these two conjectures, Moore and Seiberg  \cite{MS1}
\cite{MS2} derived a set of polynomial equations. The Verlinde conjecture and the Verlinde
formula \cite{V} were shown in \cite{MS1}
to be consequences of this set of equations but, since these two 
conjectures above were not proved in \cite{MS1}
\cite{MS2},  they remained to be conjectures until they were proved 
by the author in \cite{H-Verlinde} in 2004. 
Moreover, Moore and Seiberg noticed in \cite{MS2} that their set of polynomial equations 
has properties analogous to those of tensor categories and 
there are also additional properties that are not satisfied by the usual 
examples of tensor categories. Later, Turaev \cite{T1} \cite{T2} formulated the
precise notion of modular tensor category. 
From the set of Moore-Seiberg polynomial equations, it is not difficult to 
obtain an abstract modular tensor category in the sense of Turaev. But one cannot 
construct a modular tensor category structure on the category of modules from 
this set of Moore-Seiberg equations or this abstract modular tensor category alone. 
Since this set of Moore-Seiberg polynomial equations are consequences of the 
two conjectures on the operator product expansion and modular invariance 
of chiral vertex operators (intertwining operators) mentioned above,  
we see that these two conjectures led to a conjecture that the category of modules for 
a vertex operator algebra for a rational conformal field theory has a natural structure of 
modular tensor category. 
The vertex operator algebra associated to the category 
$\widetilde{\mathcal{O}}_{\ell}$ for a finite-dimensional complex simple Lie algebra
 $\mathfrak{g}$ and $\ell\in \Z_{+}$ was conjectured to be such a vertex operator algebra. 

These two conjectures on the operator product expansion and modular invariance 
of intertwining operators were proved by the author in 2002 \cite{H-diff-eqn} and 2003
\cite{H-modular}, respectively, for vertex operator algebras satisfying suitable 
positive-energy, finiteness and reductivity conditions.

\section{Tensor category structures on suitable categories of modules 
for affine Lie algebras}

In this section, we review the constructions of tensor category structures 
on suitable categories of modules for affine Lie algebras.

\subsection{Tensor categories}

Before we discuss the results on the constructions of tensor 
category structures on suitable categories of modules for affine Lie algebras,
we recall briefly here some basic terms we use 
in the theory of tensor categories to avoid confusions. 
See \cite{T2} and \cite{BK} for details. For definitions and terms in 
the theory of vertex tensor categories, we refer the reader to  \cite{HL2} and \cite{H-rigidity}.

A tensor category is an abelian category with a tensor product bifunctor,
a unit object, an associativity isomorphism, a left unit isomorphism and a right 
unit isomorphism such that the pentagon and triangle diagram are commutative.

A braided tensor category is a tensor category with a braiding isomorphism
such that two hexagon diagrams are commutative. 

A braided tensor category is rigid if every object has a  two-sided dual object. 

A ribbon tensor category is a rigid braided tensor category with a twist 
satisfying the balancing axioms.

A modular tensor category is a semisimple ribbon tensor category with finitely many 
simple (irreducible) objects such that 
the matrix of the Hopf link invariants 
is invertible.

\subsection{The case of $\ell+h^{\vee}\not\in \Q_{\ge 0}$}

In this case, we have the following major result of Kazhdan and Lusztig \cite{KL1}--\cite{KL5}:

\begin{thm}[Kazhdan-Lusztig]
Let $\ell\in \C$ such that $\ell+h^{\vee}\not\in \Q_{\ge 0}$. 
Then $\mathcal{O}_{\ell}$ has a natural rigid braided tensor 
category structure. Moreover, this rigid braided tensor 
category is equivalent to the rigid braided tensor 
category of finite-dimensional integrable modules for a quantum group 
constructed from $\mathfrak{g}$
at $q=e^{\frac{i\pi}{\ell+h^{\vee}}}$. 
\end{thm}

This result was announced in 1991 \cite{KL1} and 
the detailed constructions  and proofs were published in 1993 \cite{KL2} \cite{KL3}
and 1994 \cite{KL4} \cite{KL5}. The construction of the rigid braided tensor 
category structure, especially the rigidity, in this work depends heavily on 
the results on the quantum group side. In particular, 
this construction cannot be adapted directly to give constructions for the module categories
at other levels. 

In 2008, using the logarithmic generalization \cite{HLZ1}--\cite{HLZ9} 
by Lepowsky, Zhang and the author of the semisimple 
tensor category theory of Lepowsky and the author \cite{HL1}--\cite{HL6}
and of the author \cite{tensor4} \cite{H-diff-eqn}, Zhang \cite{Z}
gave a vertex-operator-algebraic  construction of the braided tensor category 
structure in this case, with a mistake corrected by the author in 2017 \cite{H-Z-correction}.
The braided tensor category obtained in this construction can be shown easily
to be indeed equivalent to the one constructed in \cite{KL1}--\cite{KL5} (for example, the 
tensor product bifunctors are in fact the same in these two constructions).
In this vertex-operator-algebraic construction, the main work is the proof of the 
associativity of logarithmic intertwining operators (logarithmic operator product expansion) 
and the construction of the associativity isomorphism. By using the work \cite{HLZ1}--\cite{HLZ9},
this work in fact constructed a vertex tensor category structure and the braided tensor category structure
was obtained naturally from this vertex tensor category structure.  But in \cite{Z}, only 
the braided tensor category structure was obtained and the rigidity is not proved. 
Thus we still have the following interesting problem:

\begin{prob}
Give a proof of the rigidity in this case in the framework of the logarithmic tensor category
theory of Lepowsky, Zhang and the author \cite{HLZ1}--\cite{HLZ9}, 
without using the results for modules for the corresponding quantum group. 
\end{prob}

\subsection{The case of $\ell\in \Z_{+}$}

 This $\ell\in \Z_{+}$ case is what the original work and conjectures of Moore and Seiberg 
\cite{MS1} \cite{MS2} were about. The construction of modular tensor category structures
in this case has a quite interesting and complicated story.

Recall from Section 3
that the work of Moore and Seiberg \cite{MS1} \cite{MS2} led to a conjecture 
that the category $\widetilde{O}_{\ell}$ has a natural structure of modular tensor 
category in the sense of Turaev \cite{T1} \cite{T2}. 

In 1997, Lepowsky and the author \cite{HL-affine}
gave a construction of the braided tensor category 
structure on the category $\widetilde{O}_{\ell}$ when $\ell\in \Z_{+}$, using the semisimple 
tensor product bifunctor constructed by Lepowsky and the author \cite{HL1}--\cite{HL6}
 and the associativity isomorphism constructed by author \cite{tensor4}
in the general setting of the category of modules for a vertex operator algebra satisfying 
suitable conditions. This work obtained the braided tensor category structure naturally from
a vertex tensor category structure constructed in  \cite{HL-affine} using 
\cite{HL3}--\cite{HL6} and \cite{tensor4}. 

In 2001, using the method developed in an early unpublished 
work \cite{BFM} of Beilinson-Feigin-Mazur in 1991,
Bakalov and Kirillov Jr. also gave in a book \cite{BK} 
a construction of the braided tensor category 
structure on the category $\widetilde{O}_{\ell}$ in this case. 
Indeed Beilinson and other people 
knew how to construct the braided tensor category structure in mid 1990's. 
But Beilinson informed the author in private discussions in 1996 that he did not know how to 
prove the rigidity. Bakalov and Kirillov Jr. also did not give a proof of 
the rigidity and the nondegeneracy 
property in \cite{BK}, though these are stated in \cite{BK} as parts of the main theorem
(Theorem 7.0.1 in \cite{BK}). 
.

Since there are a lot of confusions, we quote the sentence about rigidity 
in \cite{BK} after the proof of Corollary 7.9.3 (which is a precise statement of 
Theorem 7.0.1 in the book): ``As a matter of fact, we have 
not yet proved the rigidity (recall that modular functor only defines weak rigidity); however, it
can be shown that this category is indeed rigid.'' This is the only place in the whole book
where the rigidity of this tensor category is explicitly mentioned. There is also no reference given
for a proof or even a sketch of a proof of 
this rigidity. In 2012, Bakalov and Kirillov Jr. informed the author 
explicitly in private communications that they do not have a proof of the rigidity.

It was until 2005 (in fact the author already gave a talk on the proof in 2004 in the 
the Erwin Schr\"{o}dinger Institute but
the paper was posted in the arXiv in 2005), the rigidity and the nondegeneracy property
were finally proved by the author in \cite{H-rigidity} as a special case of a proof
for a vertex operator
algebra satisfying natural existence-of-nondegenerate-invariant-bilinear form, 
positive-energy, finiteess and reductivity conditions. In particular, the conjecture mentioned above 
was proved by the author  in 2005: 

\begin{thm}[\cite{H-rigidity}]
Let $\ell\in \Z_{+}$. Then the category $\widetilde{\mathcal{O}}_{\ell}$ 
has a natural structure of a modular tensor category.
\end{thm}

The author's proof of this theorem was based on a formula used 
by author \cite{H-Verlinde} to derive the Verlinde formula. 
This was the first indication that the rigidity is in fact deeply related to the 
Verlinde formula which in turn
is a consequence of the operator product expansion (associaitivity)  
and modular invariance of intertwining operators.

Assuming the existence of the rigid braided tensor category structure 
on $\widetilde{\mathcal{O}}_{\ell}$, 
Finkelberg in his thesis \cite{F1} in 1993 and then in a revision \cite{F2}
in 1996 gave a proof that this rigid braided tensor category
is equivalent to a semisimple subquotient of a rigid braided tensor category of modules for 
a quantum group constructed from $\mathfrak{g}$, using the equivalence constructed 
by Kazhdan and Lusztig \cite{KL1}--\cite{KL5}. This work had been reinterpreted as 
giving a construction of the rigid braided tensor category structure
on $\widetilde{\mathcal{O}}_{\ell}$. 

But in 2012,  the author found a gap in Finkelberg's paper. Through Ostrik, the author 
informed Finkelberg this gap and explained that in \cite{H-Verlinde} 
the Verlinde formula is needed to prove the rigidity. Then instead of obtaining
the Verlinde formula in the affine Lie algebra case as a major consequence,
Finkelberg filled the gap in \cite{F3} by using this formula whose proofs had been given 
using other methods by Faltings \cite{Fa}, Teleman \cite{Te}
and the author \cite{H-Verlinde}. 

Even after the correction in 2013, Finkelberg's proof of the rigidity is not complete. 
There are a few cases, including the important $\mathfrak{g}=\mathfrak{e}_{8}$ and $\ell=2$ case, 
that his method does not work. Finkelberg's equivalence between 
the modular tensor category $\widetilde{\mathcal{O}}_{\ell}$ and a semisimple 
subquotient of a rigid braided tensor category of modules for a quantum group is also not complete
because of the same few cases, including the $\mathfrak{g}=\mathfrak{e}_{8}$ and $\ell=2$ case,
in which his method does not work. Thus we have the following 
open problem:

\begin{prob}
Find a direct construction of this equivalence without using the equivalence given by 
Kazhdan-Lusztig so that this equivalence covers all the cases, including the important 
$\mathfrak{g}=\mathfrak{e}_{8}$ and $\ell=2$ case.
\end{prob}

\subsection{The admissible case}

In the case that $\ell$ is an admissible number, Creutzig, Yang and the author proved
in 2017 the following result:

\begin{thm}[\cite{CHY}]\label{CHY1}
Let $\ell$ be an admissible number. Then the category $\mathcal{O}_{\ell, {\rm \scriptsize ord}}$ 
has a natural structure of a braided tensor category with a twist.
\end{thm}

This theorem was proved using the logarithmic tensor category theory 
of Lepowsky, Zhang and the author \cite{HLZ2}--\cite{HLZ9}, some results
of Kazhdan-Lusztig \cite{KL2} and some results of Arakawa \cite{A}. 
The logarithmic tensor category theory reduces the construction
of such a braided tensor category structure to the verification of several conditions. 
Since it used \cite{HLZ2}--\cite{HLZ9} to obtain the braided tensor category structure,
this work also gave a construction of a vertex tensor category structure. 

This is a semisimple category. At first one might want to use
the early tensor category theory for semisimple category of modules 
by Lepowsky and the author \cite{HL1}--\cite{HL6} and by the author \cite{tensor4}
and \cite{H-diff-eqn}. The main result in \cite{tensor4}
in this semisimple theory constructing the associativity isomorphism
in this theory needs a convergence and extension property 
without logarithm. If generalized modules (not necessarily lower bounded)
for the affine Lie algebra 
vertex operator algebras in this case are
all complete reducible, a result in \cite{H-diff-eqn} can be applied to this case
to conclude that the convergence and extension property
without logarithm holds. But in this case, we do not have such a 
strong complete reducibility theorem
and thus we cannot directly use this semisimple theory.

Instead, we use the logarithmic generalization of the semisimple 
theory, even though our theory is semisimple.
In this theory, we do not need to prove that there is no logarithm in the analytic
extension of 
the product of intertwining operators. 
We construct the associativity isomorphism from logarithmic intertwining operators. 
Finally, since the modules in the category are all semisimple, 
the logarithmic intertwining operators are all ordinary. In particular, our theory 
still has no logarithm.

Another condition that needs to be verified is that the category should be closed under
a suitable tensor product operation. This condition is verified using a result of Arakawa in \cite{A}.

The most subtle condition is the condition that suitable submodules in the dual space of 
the tensor product of two modules in $\mathcal{O}_{\ell, {\rm \scriptsize ord}}$ should also be in 
$\mathcal{O}_{\ell, {\rm \scriptsize ord}}$. 
The verification of this condition uses the author's modification 
in 2017 of one main result in \cite{HLZ8} that had been used to correct a mistake in Zhang's construction \cite{Z} in the case 
of $\ell+h^{\vee}\not\in \Q_{\ge 0}$. It also uses some results of Kazhdan-Lusztig \cite{KL2} and Arakawa \cite{A}.

Theorem \ref{CHY1} gives only a braided tensor category structure with a twist 
to the semisimple category $\mathcal{O}_{\ell, {\rm \scriptsize ord}}$. The 
natural question is whether
this category is rigid and modular. In the case of $\mathfrak{g}=\mathfrak{sl}_{2}$,
Creutzig, Yang and the author gave a complete answer in the same paper:

\begin{thm}[\cite{CHY}]\label{CHY2}
Let $\mathfrak{g}=\mathfrak{sl}_{2}$ and $\ell=-2+\frac{p}{q}$ with $p, q$ coprime positive integers.
Then the braided tensor category $\mathcal O_{\ell,{\rm \scriptsize ord}}$ is a 
ribbon tensor category and 
is a modular tensor category if and only if $q$ is odd.  
\end{thm}

The idea of the proof of this theorem
is roughly speaking the following:
First use the theory of vertex operator algebra extensions established by Kirillov Jr., Lepowsky
and the author \cite{HKL} and by Creutzig, Kanade and McRae \cite{CKM}  to
show that this tensor category is  braided equivalent to a full tensor subcategory 
of the braided tensor category of modules for a minimal Virasoro vertex operator algebra. 
But the braided tensor category of modules for the 
minimal Virasoro vertex operator algebra in fact has a natural modular tensor category structure
constructed by the author
\cite{H-rigidity} as a special case of a general construction for 
vertex operator algebras corresponding to rational conformal field theories. 
Then Theorem \ref{CHY2} is proved by 
using this modular tensor category structure and the braided 
equivalence mentioned above. 

For general complex semisimple $\mathfrak{g}$, conjectures have been given 
in the paper \cite{CHY}. 
Let  $\ell$ an admissible number. 
Then in particular, $\ell=-h^{\vee}+\frac{p}{q}$ with coprime positive integers $p, q$. 

\begin{conj}\label{CHY-conj-1}
The braided tensor category structure on $\mathcal{O}_{\ell, {\rm \scriptsize ord}}$ 
is rigid and thus is a ribbon tensor category. 
\end{conj}

In the case that $\mathfrak{g}$ is simple-laced, Conjecture \ref{CHY-conj-1} has been proved
recently by Creutzig:

\begin{thm}[\cite{C}]
When $\mathfrak{g}$ is simple-laced, the braided tensor category structure 
on $\mathcal{O}_{\ell, {\rm \scriptsize ord}}$ 
is rigid and thus is a ribbon tensor category. 
\end{thm}

The idea of the proof of this result is similar to the proof 
in the case of $\mathfrak{sl}_{2}$ discussed above but is much harder because 
it uses the ribbon tensor category (actually a modular tensor category) 
of modules for a vertex operator algebra called the rational principle 
$W$-algebra of $\mathfrak{g}$, instead of the ribbon tensor category of 
modules for a minimal Virasoro vertex operator algebra.
Arakawa proved in \cite{A2} 
that such a $W$-algebra satisfies the conditions needed in Theorem 4.6 in \cite{H-rigidity}
and thus by Theorem 4.6 in \cite{H-rigidity},
the category of modules for this $W$-algebra has a natural 
structure of a modular tensor category. Using 
the work \cite{ACL} of Arakawa, Creutzig and Linshaw on principal $W$-algebras 
of simply-laced Lie algebras, the theory of vertex operator algebra extensions 
in \cite{HKL} and \cite{CKM} mentioned above and suitable fusion rule results
proved by Creutzig in \cite{C}, 
Creutzig constructed in \cite{C} a fully faithful braided tensor 
functor from a subcategory of the modular tensor category 
of modules for this $W$-algebra to 
a certain braided tensor category with a twist. 
The rigidity of the tensor category 
of modules for the $W$-algebra gives the rigidity of this braided tensor category 
with a twist so that it becomes a ribbon tensor category.
Since $\mathcal{O}_{\ell, {\rm \scriptsize ord}}$ is different from 
this ribbon tensor category by just the actions of certain simple currents,
$\mathcal{O}_{\ell, {\rm \scriptsize ord}}$ is also rigid and thus a ribbon tensor category.

\begin{conj}\label{CHY-conj-2}
The ribbon tensor category structure on $\mathcal{O}_{\ell, {\rm \scriptsize ord}}$ is modular
except for the following list:
\begin{enumerate}

\item $\mathfrak g \in \{\mathfrak{sl}_{2n}, \mathfrak{so}_{2n}, \mathfrak{e}_7, \mathfrak{sp}_n\}$ and $q$ even.

\item $\mathfrak g = \mathfrak{so}_{4n+1}$ and $q=0 \mod 4$.

\item $\mathfrak g = \mathfrak{so}_{4n+3}$ and $q=2 \mod 4$. 
\end{enumerate}
\end{conj}

These conjectures follow from the following
conjecture on the equivalence of these braided tensor categories 
with the braided tensor categories coming from module categories for quantum groups constructed from 
the same finite-dimensional simple Lie algebra $\mathfrak{g}$: 

\begin{conj}\label{CHY-conj-3}
The category $\mathcal O_{\ell, {\rm \scriptsize ord}}$ and 
the semi-simplification $\mathcal C_\ell(\mathfrak g)$ of the category of tilting modules 
for $U_q(\mathfrak g)$ are equivalent as braided tensor categories, where 
$q=e^{\frac{\pi i }{r^{\scriptsize \vee}(\ell+h^{\scriptsize \vee})}}$. 
\end{conj}

Since the rigidity and modularity of 
the braided tensor category of tilting modules 
for $U_q(\mathfrak g)$ were established by Sawin in 2003 \cite{S},   we see that 
Conjectures \ref{CHY-conj-1} and 
\ref{CHY-conj-2} are indeed consequences of Conjecture \ref{CHY-conj-3}.

The results and conjectures in \cite{CHY} are all about the subcategory 
$\mathcal O_{\ell, {\rm \scriptsize ord}}$ of $\mathcal{O}_{\ell}$. We still have the
following open problem:

\begin{prob}
Let $\ell$ be an admissible number. 
Is there a braided tensor category structure on $\mathcal{O}_{\ell}$? If there is, 
is it rigid or even modular in a suitable sense? 
\end{prob}

\subsection{The remaining case: an open problem}

For the remaining case, we have only an open problem:

\begin{prob}
Let $\ell$ be a rational number larger than or equal to $-h^{\scriptsize \vee}$
but is not admissible. What tensor category structures can we construct?
Are they rigid, semisimple or modular?
\end{prob}

\bigskip

\noindent {\small \sc Department of Mathematics, Rutgers University,
Piscataway, NJ 08854}

\noindent {\em E-mail address}: yzhuang@math.rutgers.edu

\end{document}